\documentclass[10pt,reqno]{amsart}
\usepackage{amssymb}
\usepackage{amsfonts}
\usepackage{amsthm}
\usepackage{amsmath}
\usepackage{epsfig}

\newcommand{\IR}{{\mathbb{R}}}
\newcommand{\IC}{{\mathbb{C}}}
\newcommand{\ZZ}{{\mathbb{Z}}}

\let\det=\undefined\DeclareMathOperator*{\det}{det}

\DeclareMathOperator*{\tr}{tr} \DeclareMathOperator*{\diag}{diag}

\theoremstyle{definition}

\newcounter{smalllist}

\begin{document}

\title{RIEMANN-HILBERT METHODS IN THE THEORY OF ORTHOGONAL POLYNOMIALS}
\author{Percy Deift}
\address{Percy Deift\\
         Courant Institute of Mathematical Sciences\\
         New York University\\
         251 Mercer Street\\
         New York, NY 10012}
\email{deift@cims.nyu.edu}


\dedicatory{To Barry Simon, on his 60th birthday.\\
Mathematician extraordinaire, teacher and friend.}

\thanks{The work of the author was supported in part by DMS Grants
No. 0500923  and No. 0296084, and also by a Friends of the Institute Visiting Membership at
the Institute for Advanced Study in Princeton, Spring 2006. The author would also
like to thank Irina Nenciu for many comments and much help.}

\begin{abstract}
In this paper we describe various applications of the Riemann-Hilbert method
to the theory of orthogonal polynomials on the line and on the circle.
\end{abstract}

\maketitle

\section{Introduction}
In this paper $d\mu$ denotes either a Borel measure on $\IR$
with finite moments
\begin{equation}\label{2.1}
\int_{\IR} |x|^m\, d \mu (x) < \infty \ , \ \ \ \ \ m = 0,1,2,\ldots
\end{equation}
or a finite Borel measure on the unit circle $S^1$
\begin{equation}\label{2.2}
\int_{S^1} d\mu (\theta) < \infty.
\end{equation}
In addition, unless stated explicitly otherwise, we will always
assume that $d\mu $ is  a {\it nontrivial probability measure}, i.e. supp$(d\mu )$ is
infinite and the integral of $d\mu$ is 1.

Let
\begin{equation}\label{2.3}
p_n (x) = k_n x^n + \ldots \ , \ \ \ \ k_n > 0 \ , \ \ \ \ n =
0,1,2,\ldots
\end{equation}
\begin{equation}\label{2.4}
\phi_n (z) = \kappa_n z^n + \ldots \ , \ \ \ \ \kappa_n > 0 \ , \ \ \ \ n =
0,1,2,\ldots
\end{equation}
denote the orthonormal polynomials (OP's) with respect to $d\mu$ on
$\IR$ and $S^1$ respectively (see \cite{Szego}),
\begin{equation}\label{2.5}
\int_{\IR} p_n (x) p_m (x) d\mu (x) = \int_{S^1}
\overline{\phi_n (e^{i\theta}) } \ \phi_m (e^{i\theta} ) d\mu
(\theta ) = \delta_{n,m} \ , \ \ \ \ \ n,m \geq 0.
\end{equation}
The fact that $d\mu$ is nontrivial implies, in particular, that the
$p_n$'s, and the $\phi_n$'s, exist and are unique for all $n \geq
0$.

As is well known, the $p_n$'s satisfy a three-term recurrence
relation
\begin{equation}\label{3.1}
b_{n-1} p_{n-1} (x) + a_n p_n (x) + b_n p_{n+1} (x) = x p_n (x) \ ,
\ \ \ \ n \geq 0
\end{equation}
where
\begin{equation}\label{3.2}
a_n \in {\IR} \ , \ \ \ \ b_n > 0 \ , \ \ \ \ n \geq 0
\end{equation}
and $b_{-1} \equiv 0$. Similarly the $\phi_n$'s satisfy the
Szeg\"{o} recurrence relation
\begin{equation}\label{3.3}
\sqrt{1 - |\alpha_n |^2 } \ \phi_{n+1} (z) = z\phi_n (z) -
\bar{\alpha}_n \phi_n^* (z) \ , \ \ \ \ n \geq 0
\end{equation}
where
\begin{equation}\label{3.4}
\alpha_n \in {\IC} \ , \ \ \ \ \ |\alpha_n | < 1 \ , \ \ \ \ \ n
\geq 0
\end{equation}
and for any polynomial $q(z) $ of degree $n$
\begin{equation}\label{3.5}
q^* (z) \equiv z^n \ \overline{q(1/\bar{z}) }
\end{equation}
denotes the so-called {\it reverse polynomial}. Following \cite{Simon 1},
we call the $\alpha_n$'s {\it Verblunsky coefficients}. A simple
computation shows that
\begin{equation}\label{3.6}
\alpha_n = - {1\over \kappa_{n+1}} \overline{\phi_{n+1} (0)} \ , \ \
\ \ \ n \geq 0 \ .
\end{equation}

On $\IR$ we define the $(n+1)\times(n+1)$ {\it Hankel determinant}
\begin{equation}\label{4.1}
D_n = {\rm det} \left( \int_{\IR} x^{j+k} d\mu (x) \right)_{0\leq
j, k\leq n} \ , \ \ \ \ \ n \geq 0 \ ,
\end{equation}
and on $S^1$ we similarly define the $(n+1)\times(n+1)$ {\it
Toeplitz determinant}
\begin{equation}\label{4.2}
\Delta_n = {\rm det} \left( \int_{S^1} e^{-i(j-k)\theta} \,d\mu
(\theta ) \right)_{0\leq j, k \leq n} \ , \  \ \ \ \ n \geq 0.
\end{equation}
The determinants $D_n$ and $\Delta_n$ are closely related to the
OP's $\{ p_n \} , \{ \phi_n \}$ respectively: Indeed one has (see for
example \cite{Szego})
\begin{equation}\label{4.3}
{D_{n-1} \over D_n} = k_n^2 \ , \ \ \ \ \ {\Delta_{n-1} \over
\Delta_n} = \kappa_n^2 \ , \ \ \ \ \ n \geq 1 \ .
\end{equation}
Given $d\mu$, the study of the algebraic and asymptotic properties
of the quantities
$$
a_n ,b_n ,p_n (x) ,k_n ,\alpha_n ,\phi_n (z) , \kappa_n \ ,
$$
and also
$$
D_n \ \ \ {\rm and} \ \ \ \Delta_n \ ,
$$
constitutes the core of the classical theory of orthogonal
polynomials.

The three-term relation \eqref{3.1} can be re-written in the form
\begin{equation}\label{5.1}
Lp(z) = zp(z) \ , \ \ \ \ p(z) = (p_0 (z) ,p_1 (z) ,p_2 (z) ,\ldots
)^T \ ,
\end{equation}
where $L$ is an infinite Jacobi matrix, i.e. $L$ is symmetric
and tridiagonal
\begin{equation}\label{5.2}
L = \left( \begin{matrix}
a_0 & b_0 & & \\
b_0 & a_1 & b_1 & 0 \\
 & b_1 & a_2 & \ddots\\
0 &  & \ddots & \ddots \end{matrix} \right)
\end{equation}
with $b_i > 0 $, $i \geq 0$. In the case that $d\mu$ has compact
support on $\IR$, the operator $L$ is bounded on
$$
\ell_2^+ = \left\{ u = (u_0 ,u_1 ,\ldots )^T : \ \sum_{i=0}^{\infty}
|u_i |^2 < \infty \right\}.
$$
Let
$$
F : \{ d\mu \ {\rm on} \ {\IR} : \ {\rm supp} (d\mu) \ {\rm
compact} \} \rightarrow \{ {\rm bounded \ Jacobi \ matrices \ on} \
\ell_2^+ \}.
$$
denote the map taking $d\mu \mapsto L$. Conversely, if $L$ is a
bounded Jacobi matrix, then in particular $L$ is self-adjoint, and
we let $d\mu$ denote the spectral measure associated with $L$ in the
cyclic subspace generated by $L$ and $e_0$, where
$e_0 = (1,0,0,\ldots )^T \in \ell_2^+$. Thus
\begin{equation}\label{6.1}
\left( e_0 ,{1\over L-\lambda } e_0 \right) = \int {d\mu (x) \over
x-\lambda} \ , \ \ \ \ \ \lambda \in {\IC} \backslash {\IR}
\end{equation}
and it follows further that $d\mu$ has compact support. Let
$$
\hat{F} : \ \{ {\rm bounded \ Jacobi \ matrices \ on} \ \ell_2^+ \}
\rightarrow \{ d\mu \ {\rm on} \ {\IR} : \ {\rm supp} (d\mu) \ {\rm
compact} \}
$$
denote the map taking $L$ to $d\mu$. The basic fact of the matter
(see, for example, \cite{Akhiezer}, \cite{Simon 2}, and also \cite{Deift 1}) is
that $F$ and $\hat{F}$ are inverse to each other, $F \circ \hat{F} =
id$, $\hat{F} \circ F = id$. From this point of view the (classical)
orthogonal polynomial problem is the inverse spectral component of a
spectral/inverse spectral problem. If the support of $d\mu$ is not
compact, then the situation is similar, but the relation between
$d\mu$ and $L$ is more complicated because $L$ is now an unbounded
operator and we must distinguish between different self-adjoint
extensions of $L$ (see \cite{Akhiezer}, \cite{Simon 2} for more details).

In the case of measures $d\mu$ on the unit circle, the role of the
Jacobi matrices is played by so-called CMV matrices $C$ (see \cite{Simon 1}).
Such matrices $C$ are unitary in $\ell_2^+$ and pentadiagonal,
and have the form
\begin{equation}\label{7.1}
C = LM
\end{equation}
where $L$ and $M$ are block diagonal
\begin{equation}\label{8.1}
L = \text{diag}(\Theta_0,\Theta_2,\Theta_4,\dots) \ , \ \ \ \
M =\text{diag}(1,\Theta_1,\Theta_3,\dots)
\end{equation}
with
\begin{equation}\label{8.2}
\Theta_j = \left( \begin{array}{cc}
\bar\alpha_j & \rho_j \\
\rho_j & -\alpha_j \end{array} \right) \  , \ \ \ \ \ j \geq 0 \ .
\end{equation}
Here
\begin{equation}\label{8.3}
|\alpha_j | < 1 \ , \ \ \ \ \ j \geq 0
\end{equation}
and
\begin{equation}\label{8.4}
\rho_j = \sqrt{ 1 - |\alpha_j |^2 }\,.
\end{equation}
CMV matrices are named for Cantero, Moral and Vel\'azquez \cite{CMV}, but
in fact they appeared earlier in the literature (see, in particular,
\cite{Watkins}). Let
$$
\psi : \ \{ d\mu \ {\rm on} \ S^1 \} \rightarrow \{ {\rm CMV
\ matrices} \}
$$
denote the map taking $d\mu \rightarrow C$, the CMV matrix
constructed from the Verblunsky coefficients  $\alpha_j = \alpha _j
(d\mu )$, $j \geq 0$, of $d\mu$, according to \eqref{7.1}, \eqref{8.1} and
\eqref{8.2}. Conversely, given a CMV matrix $C$, let $d\mu$ be the
spectral measure associated with $C$ in the cyclic subspace
generated by $C$, $C^*=C^{-1}$ and $e_0$.
Let
$$
\hat{\psi} : \ \{ {\rm CMV \ matrices} \} \rightarrow \{ d\mu \ {\rm
on} \ S^1 \}
$$
denote the map taking $C$ to $d\mu$. Then, as above (see \cite{Simon 1}),
$\psi$ and $\hat{\psi}$ are inverse to each other, and we see again
that the classical orthogonal polynomial problem on $S^1$ is
the inverse spectral component of a spectral/inverse spectral
problem.

The techniques used to analyze the direct spectral maps, $\hat{F}$
and $\hat{\psi}$, are generally very different from the techniques
used to analyze the inverse spectral maps, $F$ or $\psi$, though
sometimes there is some overlap (see e.g. \cite{Deift-Killip}). It is
also interesting to note that in the solution of integrable systems
one needs knowledge of {\it both} $\hat{F}$ and $F$ (or $\hat{\psi}$
and $\psi$). For example, the Toda lattice induces a flow $L_0
\mapsto L = L(t)$ on Jacobi matrices (\cite{Flaschka})
\begin{equation}\label{10+.1}
\begin{aligned}
{dL\over dt}
&= B(L) L - LB(L)\\
&L(t = 0) = L_0
\end{aligned}
\end{equation}
where
$$
L = \left( \begin{array}{cccc}
a_0 & b_0 & & 0 \\
b_0 & a_1 & b_1 & \\
 & b_1 & a_2 & \ddots\\
0 &  & \ddots & \ddots \end{array} \right) \ , \ \ \ \ B(L) =
\left( \begin{array}{cccc}
0 & b_0 & & 0 \\
-b_0 & 0 & b_1 & \\
 & -b_1 & 0 & \ddots\\
0 &  & \ddots & \ddots \end{array} \right) \ ,
$$
and the solution of \eqref{10+.1} is given by the following well-known
procedure (\cite{Moser}):
$$
L_0 \stackrel{\hat{F}}{\rightarrow} d\mu_0 = \hat{F} (L_0 )
\rightarrow d\mu_t (\lambda ) = {e^{2\lambda t} d\mu_0 (\lambda)
\over \int_{\IR} e^{2x t} d\mu_0 (x)}
\stackrel{F}{\rightarrow} L(t) = F(d\mu_t )
$$
The analysis of $\hat{F}$ and $\hat\psi$ has benefited greatly from the
powerful developments that have taken place over many years in the spectral theory
of Schr\"odinger operators and their discrete analogs, reaching, over the
last 20 years or so, and in the case of one dimension, a state of great precision.
Here Barry Simon and his school have played a decisive role, and we refer
the reader to \cite{Simon 1},
in particular, Part 2. The systematic analysis of $F$ begins with
the classic memoir of Stieltjes 1894-1895. Up till that point, a
great deal of information had been obtained concerning particular
polynomials, such as Legendre polynomials, Jacobi polynomials,
Hermite polynomials, etc., but a unified point of view based on the
orthogonality relation \eqref{2.5} had not yet emerged. The analysis of
$\psi$ began in 1920, when Szeg\"{o} initiated the systematic study
of polynomials orthogonal with respect to a measure on
$S^1$, as in \eqref{2.5}. Szeg\"{o}'s work in turn has led to many
remarkable developments by researchers from all over the world,
particularly the former USSR, Europe and the USA. We refer the
reader to Simon's book \cite{Simon 1}, where these developments are
discussed in great detail together with many fascinating anecdotes
concerning their discovery. Starting in the early 1950's with the
celebrated work of Gel'fand and Levitan, various techniques were
developed to recover one-dimensional Schr\"{o}dinger operators from
their spectral measures. In the 1970's, techniques based on the
inverse-Schr\"{o}dinger method (see \cite{CGe} and \cite{Flaschka}) started to play a role in the
analysis of $F$ and $\psi$. The goal of this paper is to describe
one of these techniques, which is different from the techniques in \cite{CGe} or \cite{Flaschka}, and
which has proved extremely fruitful, {\it viz.}, the Riemann-Hilbert (RH) method, also referred
to as the Riemann-Hilbert Problem (RHP). The scope of the paper is limited to
describing results for $F$ and $\psi$ obtained by RHP. Some of the results that we describe
are quite standard and are included only for purposes of illustration. Other results,
particularly asymptotic results, have been obtained, so far, only through RH methods. For a full
up-to-date discussion of what is known about $F$ and $\psi$, including the seminal
contributions of Golinskii, Ismail, Khrushchev, Lubinsky, Nevai, Rakhmanov,
Saff, Totik and many others, we again refer the reader to
\cite{Simon 2} and \cite{Simon 1}.

\smallskip

\begin{figure}\label{contour}
\begin{center}
\epsfig{file=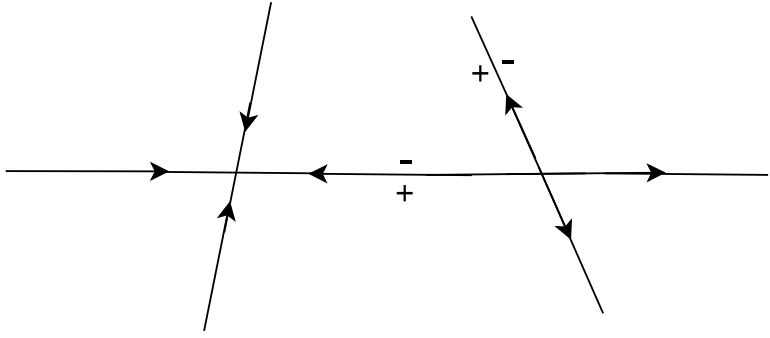,width=4in}
\end{center}
\caption{The contour $\Sigma$}
\end{figure}

To begin, let $\Sigma$ be an oriented
contour in the complex plane $\IC$ (see Figure \ref{contour}).
By convention, if we move along the contour in the direction of the
orientation, the (+)-side (resp. ($-$)-side) of the contour lies to
the left (resp. right) (see again Figure 1). A $k\times k$ {\it jump
matrix} $v$ on $\Sigma$ is a mapping from $\Sigma\rightarrow G\ell
(k,{\IC})$ such that $v,v^{-1} \in L^{\infty} (\Sigma )$. We say
that an $\ell \times k$-valued matrix function $m(z)$ is a solution
of the RHP $(\Sigma ,v)$ if
\begin{enumerate}
\item[(a)]
$m(z)$ is analytic in ${\IC} \backslash \Sigma$
\item[(b)]
$m_+ (z) = m_- (z) v(z)$, $z \in \Sigma$, where $\displaystyle{ m_{\pm}
(z) = \lim_{ z^{\prime} \rightarrow z, z \in (\pm )-{\rm side}}
m(z^{\prime}) }$
\end{enumerate}
If in addition $\ell = k$ and
\begin{enumerate}
\item[(c)]
$m (z) \rightarrow I$ as $z \rightarrow \infty$
\end{enumerate}
we say that $m$ is a solution of the {\it normalized} RHP $(\Sigma
,v)$.

\smallskip

Many technical issues arise. For example, in what sense do the
limits $m_{\pm}$ exist? In what sense does $m(z) \rightarrow I$ in
(c)? How should one understand (b) at points of self-intersection
in $\Sigma$? Under what assumptions on $\Sigma$ and $v$ does a
solution $m(z) $ exist? And if we normalize as in (c), is the
solution unique? We will not consider such issues here and in the
text that follows, and we simply refer the reader to \cite{CG} and the
references therein for a general discussion of RHP's (see also \cite{DZ1}
for more recent information, and \cite{BDT} for a discussion of points of
self-intersections). In this paper we will consider almost exclusively
problems with solutions $m(z)$ that are analytic in ${\IC}
\backslash \Sigma$ and continuous up to the boundary and at $z =
\infty$. For such solutions, the limits in (b) and (c) are taken
pointwise. Furthermore, for the problems we consider, the solution of
the normalized RHP will always exist and be unique.

At the analytical level, a normalized RHP is equivalent to a problem
for coupled singular integral equations on $\Sigma$. This is seen as follows.

Let $C^{\Sigma}$ denote the Cauchy operator on $\Sigma$,
\begin{equation}\label{16.1}
C^{\Sigma} h (z) \equiv \int_{\Sigma} {h(s) \over s-z} \ {ds\over 2\pi i}
\ , \ \ \ \ \ z \in {\IC} \backslash \Sigma
\end{equation}
with boundary values
\begin{equation}\label{16.2}
(C_{\pm}^{\Sigma} h) (z) = \lim_{\stackrel{z^{\prime} \rightarrow z}{
z^{\prime} \in (\pm )-{\rm side} } } (C^{\Sigma} h)(z^{\prime} ) \ ,
\ \ \ \ \ z \in \Sigma \ .
\end{equation}
Under reasonable conditions on $\Sigma$, $C_{\pm}^{\Sigma} \in
{\mathcal L} (L^p (\Sigma )$, the bounded operators from $L^p
(\Sigma ) \rightarrow L^p (\Sigma )$, for any $1 < p < \infty$, and
we have the relation
\begin{equation}\label{16.3}
C_+^{\Sigma} - C_-^{\Sigma} = 1 \ .
\end{equation}
Let
\begin{equation}\label{16.4}
v(z) = (v_- (z))^{-1} v_+ (z) \ , \ \ \ \ \ z \in \Sigma
\end{equation}
be any pointwise factorization of $v$ where
\begin{equation}\label{16.5}
v_{\pm} , \ (v_{\pm} )^{-1} \in L^{\infty} (\Sigma ) \ .
\end{equation}
Set
\begin{equation}\label{17.1}
\left\{ \begin{array}{c}
w_+ = v_+ - I , \ \ w_- = I-v_- \\
w  = (w_+ , w_- ) \end{array} \right.
\end{equation}
and define the singular integral operator on $\Sigma$
\begin{equation}\label{17.2}
C_w^{\Sigma} h  \equiv C_+^\Sigma (hw_- ) + C_-^{\Sigma} (hw_+ )
\end{equation}
for row $k$-vectors $h$. As $w_{\pm} \in L^{\infty} (\Sigma )$,
$C_w^{\Sigma} \in {\mathcal L} (L^p (\Sigma ))$, $1 < p < \infty$.
Suppose in addition that
\begin{equation}\label{17.3}
w_{\pm} \in L^p (\Sigma ) \ \ {\rm for \ some} \ \ \ 1 < p < \infty
\ ,
\end{equation}
and consider the equation for a $k\times k$-matrix function $\mu$
\begin{equation}\label{17.4}
(1-C_w^{\Sigma} )\mu = I
\end{equation}
in $I+L^p (\Sigma )$, or more precisely
\begin{equation}\label{17.5}
(1-C_w^{\Sigma} )\nu = C_w^{\Sigma} I = C_+^{\Sigma } w_- +
C_-^{\Sigma} w_+ \in L^p (\Sigma )
\end{equation}
where
\begin{equation}\label{17.6}
\mu = I+\nu \ , \ \ \ \ \ \nu \in L^p.
\end{equation}
If a solution $\mu = I+\nu$ of \eqref{17.4}--\eqref{17.6} exists, set
\begin{equation}\label{18.1}
m(z) = I + C^{\Sigma} (\mu (w_+ + w_- )) (z).
\end{equation}
Then a simple calculation shows that $m_{\pm} = \mu v_{\pm}$, and
hence $m_+ = m_- v$, and as $m(z) \rightarrow I$ as $z \rightarrow
\infty$, we see that \eqref{18.1} gives a solution of the normalized RHP $(\Sigma ,v)$.
Thus the normalized RHP $(\Sigma ,v )$ reduces to the analysis of the singular
integral equations \eqref{17.4}.

The connection between the OP problem and the RHP is due to Fokas,
Its and Kitaev \cite{Fokas-Its-Kitaev}. Let
\begin{equation}\label{18.2}
P_n = {1\over k_n } p_n = x^n + \ldots \ , \ \ \ \ \ n \geq 0
\end{equation}
denote the monic orthogonal polynomials associated with a measure
\begin{equation}\label{19.1}
d\mu ( x) = w(x) dx \ , \ \ \ \ \ w(x) \geq 0
\end{equation}
absolutely continuous with respect to Lebesgue measure on $\IR$,
with $x^j w(x) \in H^1 (\IR)$, the first Sobolev space,
for all $j \geq 0$. Let $\Sigma = {\IR}$, oriented from $-\infty$
to $+\infty$, and equipped with jump matrix
\begin{equation}\label{19.2}
v = v(x) = \left( \begin{array}{cc}
1 & w(x) \\
0 & 1 \end{array} \right) \ , \ \ \ \ \ -\infty < x < \infty.
\end{equation}
Finally, for any $n \geq 0$, let $X^{(n)} = (X_{ij}^{(n)}
)_{1\leq i,j\leq 2}$ solve the RHP $({\IR},v)$
$$
X^{(n)} (z) \ \ {\rm analytic \ in} \ \ {\IC} \backslash {\IR}
$$
\begin{equation}\label{19.3}
X_+^{(n)} (z) = X_-^{(n)} (z) v(z) \ , \ \ \ \ \ z \in {\IR}
\end{equation}
normalized so that
$$
X^{(n)} (z) \left( \begin{array}{cc}
z^{-n} & 0 \\
0 & z^n \end{array} \right) \rightarrow I \ \ {\rm as} \ \
z\rightarrow\infty \ .
$$
Then (\cite{Fokas-Its-Kitaev}, in addition see \cite{Deift 1}) direct
computation shows that
\begin{equation}\label{20.1}
X^{(n)} (z) = \left( \begin{matrix}
P_n (z) & C(P_n w)(z) \\
-2\pi i \ k_{n-1}^2 P_{n-1} (z) & -2\pi i \ k_{n-1}^2 C(P_{n-1}
w)(z)
\end{matrix} \right)
\end{equation}
where $C = C^{\IR} $ denotes the Cauchy operator on $\Sigma = {\bf
R}$. In particular,
\begin{equation}\label{20.2}
P_n (z) = X_{11}^{(n)} (z) \ .
\end{equation}
Furthermore, if $X_1^{(n)}$ denotes the residue of
$X^{(n)} \left( \begin{array}{cc}
z^{-n} & 0 \\
0 & z^n \end{array} \right) $ at infinity,
$$ X^{(n)}
(z)
 \left( \begin{array}{cc}
z^{-n} & 0 \\
0 & z^n \end{array} \right)  = I + {X_1^{(n)}
\over z} + O\left( {1\over z^2} \right) ,
$$
then
\begin{equation}\label{20.3}
k_{n-1}^2 = - {1\over 2\pi i} (X_1^{(n)} )_{21}
\end{equation}
and in the notation of \eqref{3.1}
\begin{equation}\label{20.4}
a_n = (X_1^{(n)} )_{11} - (X_1^{(n+1)} )_{11}
\end{equation}
\begin{equation}\label{20.5}
b_{n-1}^2 = (X_1^{(n)} )_{12}  (X_1^{(n+1)} )_{21}
\end{equation}
Also by \eqref{4.3} and \eqref{20.3},
\begin{equation}\label{20.6}
{D_{n-1} \over D_n} = - {1\over 2\pi i} (X_1^{(n+1)} )_{21}
\end{equation}
Thus all the basic quantities of interest in the OP problem can be
read off from the solution $X^{(n)}$ of the RHP $({\IR}, v)$
above.

On the unit circle, the situation is similar. Let
\begin{equation}\label{21.1}
\Phi_n = {1\over \kappa_n} \phi_n = z^n + \ldots \ , \ \ \ \ \ n
\geq 0
\end{equation}
denote the monic orthogonal polynomials associated with a measure
\begin{equation}\label{21.2}
d\mu (\theta ) = \omega(\theta ) {d\theta \over 2\pi}
\end{equation}
absolutely continuous with respect to Lebesgue measure on
$S^1$ with $\omega (\theta ) \in H^1 ( S^1 )$,
$\omega (\theta ) = \omega (\theta + 2\pi )$. Fix $n\geq 0$ and let
$\Sigma = S^1$, oriented counterclockwise. Equip $S^1$ with the
jump matrix
\begin{equation}\label{21.3}
v = v(\theta ) = \left( \begin{array}{cc}
1 & \omega (\theta )z^{-n} \\
0 & 1 \end{array} \right) \ , \ \ \ \ \ z = e^{i\theta}
\end{equation}
and let $Y^{(n)} = (Y_{ij}^{(n)} )_{1\leq i, j\leq 2}$ solve the RHP $(S^1 ,v)$
\begin{align}
&\bullet \ \ \ Y^{(n)} (z) \ \ {\rm analytic \ in} \ \ \IC\setminus S^1\\
&\bullet \ \ \ Y_+^{(n)} (z) = Y_-^{(n)} (z) v(\theta ) \ , \ \ \ z=e^{i\theta} \in S^1\label{22.0}
\end{align}
normalized so that
\begin{equation}
\bullet \ \ \ Y^{(n)} (z) \left( \begin{array}{cc}
z^{-n} & 0 \\
0 & z^n \end{array} \right) \rightarrow I \ \ {\rm as} \ \
z\rightarrow\infty \ .
\end{equation}
Then again (cf. \cite{BDJ}) direct computation shows that
\begin{equation}\label{22.1}
Y^{(n)} (z) = \left( \begin{matrix}
\Phi_n (z) & C(\Phi_n \ \omega/S^n ) (z) \\
-\kappa_{n-1}^2 \Phi_{n-1}^* (z) & -\kappa_{n-1}^2 C(\Phi_{n-1}^* \
\omega/S^n )(z) \end{matrix} \right)
\end{equation}
where $C = C^{S^1}$ denotes the Cauchy operator on
$S^1$ and $\Phi_{n-1}^*$ is the reverse polynomial as in
\eqref{3.5}. In particular
\begin{equation}\label{22.2}
\Phi_n (z) = Y_{11}^{(n)} (z)
\end{equation}
and hence by \eqref{3.6},
\begin{equation}\label{22.3}
\alpha_{n-1}  = - \overline{Y_{11}^{(n)} (z=0)} \ .
\end{equation}
Also
\begin{equation}\label{23.1}
\kappa_{n-1}^2 = - Y_{21}^{(n)} (z = 0)
\end{equation}
and hence
\begin{equation}\label{23.2}
{\Delta_{n-2} \over \Delta_{n-1} } = - Y_{21}^{(n)} (z=0).
\end{equation}
Again we see that all basic quantities in the OP problem on the
circle are expressed in terms of the solution $Y^{(n)}$ of the RHP
$(S^1,v)$.

\smallskip

The outline of the paper is as follows. In Section 2 we show how to use
the RHP's $({\IR},v)$ and
$(S^1,v)$ above to derive various identities, equations and
formulae for the OP problem. In Section 3 we describe the
application of the steepest descent method of Deift-Zhou for RHP's
to asymptotic problems for OP's. Finally, in Section 4 we describe
the application of RH ideas to areas related to the OP problem, such
as random matrix theory, multi-orthogonal polynomials, orthogonal Laurent polynomials,
and the rarefaction problem for the Toda lattice.

\smallskip

\noindent {\bf Technical Remark}
In most of the paper we will be considering probability measures
with some degree of smoothness as in \eqref{19.1} and \eqref{21.2} above. For
such weights we then use the RHP's to derive, in particular, various
identities such as \eqref{3.1}, \eqref{3.3}, \eqref{37.1}, etc. If $d\mu (x)$ is an
arbitrary probability measure on ${\IR}$ with finite moments, or
$d\mu (\theta )$ is a probability measure on $S^1$, we can
approximate $d\mu (x)$ and $d\mu (\theta )$ appropriately with
smooth measures $d\mu_{\epsilon} (x)$ and $d\mu_{\epsilon} (\theta
)$ respectively: For such measures \eqref{3.1}, \eqref{3.3}, \eqref{37.1}, etc., are
true, and letting $\epsilon \downarrow 0$ we conclude that these
identities are true, as they should be, for all measures $d\mu (x)$
and $d\mu (\theta )$ as above. Similar considerations apply at many
points in the paper and we leave the details to the interested
reader.

\section{Applications of $({\IR},v)$ and $(S^1,v)$: identities, equations and formulae}

The applications of Riemann-Hilbert techniques to OP's are
principally of two types:

\begin{enumerate}
\item[(a)] algebraic
\item[(b)] asymptotic.
\end{enumerate}

Under (a), the goal is to derive identities, equations and useful
formulae for the OP problems. Under (b), the goal is to determine
the asymptotic behavior of the OP's $p_n$, $P_n$, $\phi_n$, $\Phi_n$ as
$n\rightarrow\infty$: Here one considers the case where the weight
$\omega (x)$ is independent of $n$, as well as the case where
$\omega (x)$ depends on $n$ in a prescribed fashion (see \eqref{53.2} below).
We consider (a) in this section, and (b) in the next.

Regarding (a), there is a general methodology, which may be traced
all the way back to the original work of Gel'fand and Levitan, and
which may be stated loosely as follows: If the jump matrix for a RHP
is independent of a parameter, then differentiation with respect to
that parameter (or taking differences in the discrete case) leads to
an equation/identity.

We illustrate this methodology, which may be viewed as the analog
for RHP's of the celebrated theorem of Noether on conserved
quantities for dynamical systems, first in the case of the
defocusing Nonlinear Schr\"{o}dinger Equation (NLS). In 1975 Shabat
observed that the inverse scattering problem for the one-dimensional
Schr\"{o}dinger equation could be rephrased as a RHP. Because of the
connection between Schr\"{o}dinger operators and the Korteweg-de
Vries (KdV) equation, this meant that KdV, and by extension all
1+1-dimensional integrable systems, could be solved by a RHP. In the
case of defocusing NLS, Shabat's observation amounts to the
following (see e.g. \cite{DZ1}). Let $q(x,t)$ be the solution of NLS on
the line
\begin{equation}\label{26.1}
\left\{ \begin{array}{c}
iq_t + q_{xx}  - 2|q|^2 q = 0 \\
q(x,t = 0) = q_0 (x) \end{array} \right.
\end{equation}
where $q_0 (x) \rightarrow 0$ sufficiently rapidly as
$|x|\rightarrow\infty$. Just as KdV is associated with the
Schr\"{o}dinger operator, NLS is associated with a first order,
two-by-two scattering problem
\begin{equation}\label{26.2}
\frac{d\psi}{dx} = i \frac{z}{2} \sigma_3 \psi +
\left(\begin{matrix}
0 & iq\\
-i\bar{q} & 0
\end{matrix}\right)
\psi \ , \qquad -\infty < x <
\infty
\end{equation}
where $\displaystyle{ \sigma_3 = \left( \begin{array}{cc}
1 & 0 \\
0 & -1 \end{array} \right) }$ is the third Pauli matrix. Let $r=
r(z)$ be the reflection coefficient for \eqref{26.2} with $q = q_0$. The
map $\hat{R}:q \mapsto r$ is the analog for NLS of the OP maps
$\hat{F}$ and $\hat{\phi}$. Now, for fixed $x$ and $t$, let $m =
m(z; x,t)$ be the solution of the normalized RHP $({\IR}
,v_{x,t})$ where $\IR$ is oriented from $-\infty$ to $+\infty$ and
\begin{equation}\label{27.1}
\left\{ \begin{array}{l} v_{x,t} (z) = \left( \begin{array}{ll}
1 - |r(z)|^2 & re^{i\theta} \\
- \bar{r} e^{i\theta} & 1 \end{array} \right) \\
\theta = xz - tz^2 \end{array} \right. \ , \ \ \ \ z \in {\IR}.
\end{equation}
Let $m_1 (x,t)$ be the residue of $m$ at $z = \infty$,
$$
m(z; x,t) = I + {m_1 (x,t) \over z} + O \left( {1\over z^2 } \right)
\ .
$$
Then
\begin{equation}\label{27.2}
q(x,t) = -i (m_1 (x,t))_{12}
\end{equation}

How does one prove \eqref{27.2}? At the functorial level, $\hat{R}$ is
really a map from the category of differential operators to the
category of RHP's,
$$
L(q) \mapsto q \mapsto r \mapsto v_{x,t}
$$
where $L(q)= i\sigma_3\frac{d}{dx}+ \left(\begin{matrix}0&iq\\-i\bar q&0 \end{matrix}\right)$, and
so the key question becomes: ``How is the differential operator
encoded into the formalism of RHP's?''

To answer this question, observe that
\begin{equation}\label{28.1}
\psi = \psi (z; x,t) \equiv m(z; x,t) e^{i{\theta\over 2} \sigma_3 }
\end{equation}
solves the RHP \vspace{.2in}

$\bullet \ \ \ \ \psi (z; x,t)$ analytic on ${\IC} \backslash {\IR}$ \vspace{.2in}

$\bullet \ \ \ \ \displaystyle{ \psi_+ =\psi_- \left(
\begin{array}{ll}
1 - |r(z)|^2 & r(z) \\
- \overline{r(z)} & 1 \end{array} \right) \ , \ \ \ \ \ z \in {\IR} }$ \vspace{.2in}

\noindent where the jump matrix is now independent of $x$ and $t$.
Differentiating with respect to $x$, we obtain
$$
\psi_{x+} = \psi_{x-} \left( \begin{array}{ll}
1 - |r(z)|^2 & r(z) \\
- \overline{r(z)} & 1 \end{array} \right)
$$
from which it follows that $T \equiv \psi_x\psi^{-1}$ has no jump across
$\IR$, and hence is entire. But as $z\rightarrow\infty$,
\begin{align*}
T & = m_x m^{-1} + m {iz \over 2} \sigma_3 m^{-1} \\
 & = iz {\sigma_3 \over 2} + A + O\left( {1\over z} \right)
\end{align*}
for some constant matrix $A$. By Liouville, we must then have $\displaystyle{
T = iz {\sigma_3 \over 2} + A }$ or
\begin{equation}\label{29.1}
\psi_x = iz {\sigma_3 \over 2} \psi + A \psi
\end{equation}
Simple symmetry considerations imply that $A$ is of the form
$\displaystyle{ \left( \begin{array}{cc} 0 & q \\ \bar{q} & 0
\end{array} \right) }$, and hence we recover the differential
equation \eqref{26.2}. Differentiating the $\psi$-RHP with respect to $t$
yields similarly an equation of the form
\begin{equation}\label{29.2}
\psi_t = B\psi
\end{equation}
for some explicit matrix $B = B(z, q,q_x )$. Cross-differentiating
\eqref{29.1} and \eqref{29.2}, $(\psi_x )_t = (\psi_t )_x$, then yields the NLS
equation \eqref{26.1}. It is in this way in general that identities and
differential relationships are encoded into the RHP.

To apply the above methodology to OP's, consider the solution
$X^{(n)}$ of the RHP $({\IR}, v)$ above. Observing that
$X^{(n+1)}$ satisfies the same jump relation as $X^{(n)}$ across $\IR$, we
conclude as before that $T \equiv X^{(n+1)} (X^{(n)} )^{-1}$ is entire.
But
$$
\begin{aligned}
T &= X^{(n+1)} (z) (X^{(n)} (z) )^{-1} \\
&=\left[\left( I + {X_1^{(n+1)} \over z} + O \left( {1\over z^2}\right) \right)
z^{(n+1)\sigma_3}\right] \left[ \left( I + {X_1^{(n)}
\over z} + O \left( {1\over z^2} \right) \right) z^{n\sigma_3}\right]^{-1} \\
& = z \left( \begin{array}{cc} 1 & 0 \\ 0 & 0 \end{array} \right)
+ X_1^{(n+1)} \left( \begin{array}{cc} 1 & 0 \\ 0 & 0 \end{array}
\right) - \left( \begin{array}{cc} 1 & 0 \\ 0 & 0 \end{array}
\right) X_1^{(n)} + O \left( {1\over z} \right)
\end{aligned}
$$
and again by Liouville we conclude that
\begin{equation}\label{31.1}
X^{(n+1)} (z) = \left( z \left( \begin{array}{cc} 1 & 0 \\ 0 & 0
\end{array} \right) + X_1^{(n+1)} \left( \begin{array}{cc} 1 & 0 \\
0 & 0 \end{array} \right) - \left( \begin{array}{cc} 1 & 0 \\ 0 & 0
\end{array} \right) X_1^{(n)} \right) X^{(n)}(z)
\end{equation}
from which the three-term recurrence relation \eqref{3.1} now follows by a
simple computation. Moreover, if we replace the weight $\omega (x)$
with $\displaystyle{ \omega_t (x) = {e^{2xt} \omega (x) \over
\int_{\IR} e^{2st} \omega (s) ds } }$, then
\begin{equation}\label{31.2}
W^{(n)} (z; t) \equiv X^{(n)} (z; t) e^{(tz +g(t))\sigma_3 } \ , \ \
\ \ g(t) \equiv - {1\over 2} \log \int_{\IR} e^{2st} \omega (s) ds
\ ,
\end{equation}
solves the RHP $({\IR}, v)$ with jump matrix $\displaystyle{ v =
\left( \begin{array}{cc} 1 & \omega (x) \\ 0 & 1 \end{array} \right)
}$ independent of $t$. Differentiating with respect to $t$, we
obtain as above a differential equation for $W^{(n)}$
$$
\frac{d}{dt} W^{(n)}= ( (z+ \dot{g} ) \sigma_3 + X_1^{(n)} \sigma_3 -
\sigma_3 X_1^{(n)} ) W^{(n)}\,.
$$
Using $\Gamma$ to denote the shift operator, $\Gamma W^{(n)} =
W^{(n+1)}$, equation \eqref{31.1} takes the form
\begin{equation}\label{32.1}
\Gamma W^{(n)} = \left(  z \left( \begin{array}{cc} 1 & 0 \\ 0 & 0
\end{array} \right) + X_1^{(n+1)}
 \left( \begin{array}{cc} 1 & 0 \\ 0 & 0 \end{array} \right) -
 \left( \begin{array}{cc} 1 & 0 \\ 0 & 0 \end{array} \right)
X_1^{(n)} \right) W^{(n)}
\end{equation}
Cross-``differentiating'' \eqref{31.2} and \eqref{32.1}, $\displaystyle{ {d\over
dt} \Gamma W^{(n)} = \Gamma {dW^{(n)} \over dt} }$, one is led
immediately to the Toda flow \eqref{10+.1}.

In another direction, if $\omega (x) = e^{-V(x)}$, $V(x) = \gamma_m
x^{2m} + \ldots$, $\gamma_m > 0$, then $U^{(n)} \equiv X^{(n)}
e^{\frac12 V(x) \sigma_3}$ satisfies a jump relation across $\IR$ with jump matrix
$\displaystyle{ v = \left( \begin{array}{cc} 1 & 1 \\ 0 & 1
\end{array} \right) }$, which is independent of $z$,
and by the above general methodology this leads to a differential equation
for $U^{(n)}$ with respect to $z$, $\displaystyle{ dU^{(n)} \over dz} =
DU^{(n)}$, for some explicit $D$. Cross-``differentiation'',
$\displaystyle{ {d\over dz} \Gamma U^{(n)} = \Gamma {dU^{(n)} \over dz} }$,
then leads to so-called ``string equations'' for the recurrence
coefficients $a_n ,b_n$.

Applying the above methodology to the RHP $(S^1 ,v)$ for OP's
on the unit circle, we obtain, in particular, simple and direct
proofs of Szeg\"{o} recurrence, Geronimus' Theorem on the Schur
iterates, and the Pinter-Nevai formula (see \cite{Simon 1}, and below).
Indeed, let $Y^{(n)}$ solve the RHP $(S^1 ,v)$ above. Then
one observes that $\displaystyle{ V^{(n)} \equiv Y^{(n+1)} \left(
\begin{array}{cc} 1 & 0 \\ 0 & z \end{array} \right) }$ satisfies
the same jump relation as $Y^{(n)}$ across $S^1$,
$$
\displaystyle{ V_+^{(n)} = V_-^{(n)} \left( \begin{array}{cc} 1 &
\omega z^{-n} \\ 0 & 1 \end{array} \right) },
$$
and hence $V^{(n)}
(Y^{(n)} )^{-1}$ is entire. As before, this leads to an equation for
$V^{(n)}$ and $Y^{(n)}$, which takes the form
\begin{equation}\label{33.1}
Y^{(n+1)}
 \left( \begin{array}{cc} 1 & 0 \\ 0 & z \end{array} \right)
= V^{(n)} =
 \left( \begin{array}{cc} z+\hat{a}_n & \hat{b}_n \\
\hat{c}_n & 1 \end{array} \right) Y^{(n)}
\end{equation}
for suitable constants $\hat{a}_n ,\hat{b}_n ,\hat{c}_n$.
Furthermore (det $Y^{(n)})_+ =$ (det $Y^{(n)})_-$ det $v =$ (det
$Y^{(n)})_-$, and so det $Y^{(n)}$ is entire. But det $Y^{(n)}=$ det
$\displaystyle{ \left( Y^{(n)}
 \left( \begin{array}{cc} z^{-n} & 0 \\ 0 & z^n \end{array} \right)
\right) \rightarrow 1 }$ as $z \rightarrow\infty$, and hence det
$Y^{(n)} \equiv 1$. Taking determinants of both sides of \eqref{33.1}, we
find the relation
\begin{equation}\label{32.2}
\hat{a}_n = \hat{b}_n \hat{c}_n\,.
\end{equation}
>From the first column of \eqref{33.1} we obtain the relations
\begin{equation}\label{32.3}
\Phi_{n+1} = (z+\hat{a}_n )\Phi_n - \kappa_{n-1}^2 \hat{b}_n
\Phi_{n-1}^*
\end{equation}
\begin{equation}\label{32.4}
-\kappa_n^2 \Phi_n^* = \hat{c}_n \Phi_n - \kappa_{n-1}^2
\Phi_{n-1}^*\,.
\end{equation}
Eliminating $\Phi_{n-1}^*$, we obtain the Szeg\"{o} recurrence
relation \eqref{3.3}
\begin{equation}\label{32.5}
\Phi_{n+1} = z\Phi_n - \bar{\alpha}_n \Phi_n^*
\end{equation}
with $\phi_n$ replaced by $\Phi_n$, and with Verblunsky coefficient
\begin{equation}\label{34.1}
\alpha_n = \kappa_n^2 \overline{\hat{b}_n}\,.
\end{equation}
Letting $z\rightarrow\infty$ in \eqref{32.4} we find
\begin{equation}\label{34.2}
\hat{c}_n = \kappa_n^2 \alpha_{n-1}
\end{equation}
and hence by \eqref{32.2}
\begin{equation}\label{34.3}
\hat{a}_n = \bar{\alpha}_n \alpha_{n-1}.
\end{equation}

Now consider the second column in \eqref{33.1}. Setting
\begin{equation}\label{34.4}
r_n = C(\Phi_n \omega s^{-n} ) \ , \ \ \ \ \ t_n = C(\Phi_n^*
\omega s^{-n-1} )
\end{equation}
and using \eqref{34.1}, \eqref{34.2} and \eqref{34.3}, we obtain as in \eqref{32.3} and
\eqref{32.4}
\begin{equation}\label{34.5}
zr_{n+1} = (z+\bar{\alpha}_n \alpha_{n-1}) r_n -
\bar{\alpha}_n \left( {\kappa_{n-1} \over \kappa_n} \right)^2
t_{n-1}
\end{equation}
\begin{equation}\label{34.6}
-z\kappa_n^2 t_n = \kappa_n^2 \alpha_{n-1} r_n - \kappa_{n-1}^2
t_{n-1}\,.
\end{equation}
Eliminating $t_{n-1}$ as we eliminated $\Phi_{n-1}^*$ above, \eqref{34.5}
and \eqref{34.6} reduce to
\begin{equation}\label{35.1}
r_{n+1} = r_n - \bar{\alpha}_n t_n
\end{equation}
\begin{equation}\label{35.2}
zt_{n+1} = -\alpha_n r_n + t_n\,.
\end{equation}
Defining
\begin{equation}\label{35.3}
f_n \equiv t_n /r_n
\end{equation}
and using \eqref{35.1} and \eqref{35.2}, we obtain the recurrence relation
\begin{equation}\label{35.4}
zf_{n+1} = {f_n - \alpha_n \over 1-\bar{\alpha}_n f_n } \ , \ \ \ \
\ n \geq 0.
\end{equation}
In particular, for $z=0$, we see that
\begin{equation}\label{35.5}
\alpha_n = f_n (0)
\end{equation}
and so \eqref{35.4} can be written in the form
\begin{equation}\label{35.6}
zf_{n+1} = {f_n - f_n (0) \over 1-\overline{f_n (0)}f_n} \ , \ \ \ \
\ n \geq 0.
\end{equation}
Finally observe that
\begin{equation}\label{36.1}
f_0 (z) = {t_0 \over r_0} = { { \int_{S^1} {\pi_0^*
\over s-z} \omega {ds \over 2\pi is} } \over { \int_{S^1}
{\pi_0 \over s-z} \omega {ds\over 2\pi i} } } = {\int_{S_1} {d\mu
(\theta )\over s-z} \over \int_{S_1} s {d\mu (\theta ) \over s-z} }
\ , \ \ \ \ s = e^{i\theta}
\end{equation}
where $\displaystyle{ d\mu (\theta ) = \omega (\theta ) {d\theta
\over 2\pi } }$.

Geronimus' Theorem (see \cite{Simon 1}) states the following: Let
$$
F(z) = \int_{S^1} {s+z\over s-z} \,d\mu (\theta)
$$
be the Carath\'eodory function for $d\mu$ and let $f_\text{Schur}\equiv
{1\over z} \frac{F(z) -1}{F(z)+1}$ be the associated Schur
function. Let $(f_n )_{n\geq 0}$ solve the recurrence relation
\eqref{35.6} with $f_n |_{n=0} = f_{\rm Schur}$. Then
$$
f_n (0) = \alpha_n \ , \ \ \ \ \ n \geq 0
$$
where $\{ \alpha_n \}_{n\geq 0}$ are the Verblunsky coefficients for
$d\mu$.

However, a simple computation shows that $f_{\rm Schur}$ is
precisely $f_0$ in \eqref{36.1}: Hence, using the general methodology
for RHP's as above, we have proved Geronimus' Theorem. Moreover we have the
following formula for the Schur iterates:
\begin{equation}\label{37.1}
f_n (z) = {t_n \over r_n} = {\int_{S^1} {\Phi_n^* s^{-n}
\over s-z} d\mu (\theta ) \over \int_{S^1} {\Phi_n s^{-n+1}
\over s-z} d\mu (\theta ) } \ , \ \ \ \ \ n \geq 0
\end{equation}
which reduces simply, using (3.2.52) and (2.2.53) \cite{Simon 1}, to Golinskii's formula
(\cite{Simon 1}, Thm. 32.7).

Finally we note from \cite{Simon 1}, (1.3.79), together with the simple
identity $\int {s\over s-z} d\mu=
1+z \int {d\mu \over s-z}$, that
\begin{equation}\label{37.2}
f_n = {f_{\rm Schur} B_{n-1} - A_{n-1} \over zB_{n-1}^* - zA_{n-1}^*
f_{\rm Schur} } = {(B_{n-1} - zA_{n-1} ) \int {d\mu\over s-z} - A_n
\over (zB_{n-1}^* - A_{n-1}^* ) \int {s d\mu \over s-z} + A_{n-1}^* }
\end{equation}
where $A_{n-1}, B_{n-1}$ are the Wall polynomials. But from \eqref{37.1}, we obtain
$$
f_n (z) = { z^n \int_{S^1} ({\Phi_n^* (s) s^{-n} -\Phi_n^*
(z)z^{-n} \over s-z}) d\mu (\theta ) + \Phi_n^* (z) \int {d\mu
(\theta )\over s-z} \over
 z^n \int_{S^1} ({\Phi_n (s) s^{-n} -\Phi_n (z)z^{-n} \over
s-z}) s d\mu (\theta ) + \Phi_n (z) \int s {d\mu (\theta )\over s-z}
}\,.
$$
Comparing with \eqref{37.2} we obtain
\begin{equation}\label{38.1}
\Phi_n^* (z) = B_{n-1} - zA_{n-1}
\end{equation}
or equivalently
\begin{equation}\label{38.2}
\Phi_n (z) = z B_{n-1}^* - A_{n-1}^*
\end{equation}
which is the Pinter-Nevai formula (see \cite{Simon 1}) relating the OP's
to the Wall polynomials.

In addition to the formulae and identities obtained above for OP's
using the RHP's $({\IR},v)$ and $(S^1 ,v)$, one can, using
RHP's closely related to $({\IR},v)$ and $(S^1 ,v)$, derive
formulae for Toeplitz and Hankel determinants, or more precisely
``relative'' Toeplitz and Hankel determinants, that are particularly
useful for asymptotic analysis. The asymptotic analysis of Toeplitz
and Hankel determinants, dating back at least to the work of
Szeg\"{o} in 1915, is of considerable, and continuing, mathematical
and physical interest, and we refer the reader to \cite{Basor-Widom},
\cite{Ehrhardt} and the references therein for more information and
recent results. The ``relative'' determinant formulae are as
follows.

Let $\omega_1 (x), \omega_2 (x) \geq 0$ be two weights on $\IR$
and let $D_n (\omega_1\omega_2 )$, $D_n (\omega_2 )$ be the Hankel
determinants associated with the measures $\omega_1 (x) \omega_2 (x)
dx$ and $\omega_2 (x)dx$ respectively. (Here we do not require $\omega_1\omega_2 dx$
and $\omega_2 dx$ to be probability measures.) Then
\begin{equation}\label{39.1}
\log {D_n (\omega_1 \omega_2 ) \over D_n (\omega_2 )} =
\int_0^1 dt \int_{\IR} R_t (x) \left( {d\over dt} \log
\omega_t (x) \right) dx
\end{equation}
where $\omega_t = 1-t+ t\omega_1 (x)$, $0 \leq t \leq 1$, and $R_t$ is expressed in
terms of the solution $X_t^{(n+1)} = ((X_t^{(n+1)} )_{ij} )_{1\leq
i, j\leq 2}$ of the RHP $({\IR} ,v_t(x) )$ in \eqref{19.3} with
$$\displaystyle{ v_t (x) =
\left( \begin{array}{cc} 1 & \omega_t (x) \omega_2 (x) \\
0 & 1 \end{array} \right) },$$
as follows:
\begin{equation}\label{40.1}
R_t (x) = {1\over 2\pi i} \Bigl(\bigl(X_t^{(n+1)}\bigr)_{11} \bigl(X_t^{(n+1)}\bigr)_{21}^{\prime}
-\bigl(X_t^{(n+1)}\bigr)_{11}^{\prime} \bigl(X_t^{(n+1)}\bigr)_{21} \Bigr)\,
\omega_t \omega_2.
\end{equation}
Similarly, if $\omega_1 (\theta )$, $\omega_2 (\theta ) \geq 0$ are
two weights on $S^1$, with associated Toeplitz determinants
$\Delta_n (\omega_1 \omega_2 )$, $\Delta_n (\omega_2 )$
respectively, then
\begin{equation}\label{41.1}
\log {\Delta_n (\omega_1 \omega_2 ) \over \Delta_n (\omega_2 ) } =
\int_0^1 dt \int_{S^1} R_t (\theta ) {d\over dt} \log
\omega_t (\theta ) {d\theta \over 2\pi}
\end{equation}
where $\omega_t (\theta ) = 1-t + t\omega_1 (\theta )$, $0 \leq t
\leq 1$, and $R_t (\theta )$ is expressed in terms of the solution
$Y_t^{(n+1)} = (( Y_t^{(n+1)} )_{ij} )_{1\leq i,j\leq 2}$ of the RHP
$(S^1 ,v_t (\theta ))$ in \eqref{22.0} with
$$\displaystyle{ v_t
(\theta ) = \left( \begin{array}{cc}
1 & \omega_t (\theta )\omega_2 (\theta )z^{-(n+1)} \\
0 & 1 \end{array} \right) },\quad z = e^{i\theta},
$$
as follows:
\begin{equation}\label{41.2}
R_t (\theta ) = \Bigl(\bigl(Y_t^{(n+1)}\bigr)_{11} \bigl(Y_t^{(n+1)}\bigr)_{21}^{\prime}
-\bigl(Y_t^{(n+1)}\bigr)_{11}^{\prime} \bigl(Y_t^{(n+1)}\bigr)_{21} \Bigr)\, {\omega_t
\omega_2 \over z^n}
\end{equation}
where $\prime \equiv {d\over dz}$.

The functions $R_t (x)$, $R_t (\theta )$ have the interpretation as
1-point functions
\begin{equation}\label{41.3}
R_t (x) = (n+1) \int_{x_i \in {\IR}, 1\leq i\leq n} d\mu (x,x_1
,x_2 ,\ldots ,x_n )
\end{equation}
\begin{equation}\label{41.4}
R_t (\theta ) = (n+1) \int_{\theta_i \in S^1, 1\leq i\leq n}
d\mu (\theta,\theta_1 ,\ldots ,\theta_n )
\end{equation}
for the random particle ensembles (see \cite{Mehta}) with distributions
\begin{equation}\label{42.1}
d\mu (x_0 ,x_1 ,\ldots ,x_n ) = (1/Z_{\IR}) \prod_{0\leq j < k \leq n}
(x_i - x_j )^2 \prod_{j=0}^n (\omega_t \omega_2 ) (x_j ) dx_0 dx_1
\ldots dx_n
\end{equation}
and
\begin{equation}\label{42.2}
d\mu (\theta_0 ,\theta_1 ,\ldots ,\theta_n ) = (1/Z_{S^1})
\prod_{0 \leq i < j \leq n} | e^{i\theta_j}
-e^{i\theta_k} |^2 \prod_{j=0}^n (\omega_t \omega_2 )(\theta_j
)d\theta_0 d\theta_1 \ldots d\theta_n
\end{equation}
where $Z_{\IR}$, $Z_{S^1}$ are normalization constants.

Note that on $S^1$ we can set $\omega_2 =1$, so that
$\Delta_n (\omega_2 ) =1$ and \eqref{41.1} gives us a formula, first
derived in \cite{Deift 2}, purely for $\Delta_n (\omega_1 )$. In the
non-compact situation on $\IR$, this clearly cannot be done and we
must always work with relative determinants as in \eqref{39.1}.

Formulae \eqref{39.1},\eqref{41.1} are due to Deift \cite{Deift 3}, and may be proved
by generalizing the proof of \eqref{41.1} given in \cite{Deift 2} for the
case $\omega_2 =1$. A key ingredient in the proof is the notion of
an {\it integrable operator}: If $\Sigma$ is an oriented contour in
$\IC$, we say that an operator $K$ acting on $L^p (\Sigma )$, $1 <
p< \infty$, is {\it integrable} if it has a kernel of the form
\begin{equation}\label{43.1}
K(z,z^{\prime} ) = { \sum_{j=1}^{\ell} f_j (z)g_j (z^{\prime} )
\over z-z^{\prime} } \ , \ \ \ \ \ z,z^{\prime} \in \Sigma
\end{equation}
for some functions $f_j ,g_k \in L^{\infty} (\Sigma )$, $1\leq j,k\leq \ell$.
Special examples of integrable operators appeared in
the 1960's in the work of McCoy, Tracy and others, and elements of
the general theory were discovered by Sakhnovich in the late 60's,
but the full general theory of such operators is due to Its, Izegin,
Korepin and Slavnov \cite{IIKS} in 1990. Integrable operators have many
useful properties (see e.g. \cite{Deift 2}). In particular, if $K$ is integrable as in \eqref{43.1} above, then so is
$(1-K)^{-1} -1$,
$$
(1-K)^{-1} = 1 + \frac{\sum_{j=1}^{\ell} F_j (z) G_j (z^{\prime} )}{z-z^{\prime}}
$$
for suitable $F_j, G_k$, $1\leq j,k\leq l$.
Furthermore, quite remarkably, the functions $F = (F_1 ,\ldots ,F_l
)^T$, $G = (G_1 ,\ldots G_l )^T$ can be computed in terms of the
solution of a canonical, auxiliary RHP. Indeed, define the jump
matrix $v = I-2\pi i fg^T$ on $\Sigma$, where $f = (f_1 ,\ldots ,f_l
)^T$, $g = (g_1 ,\ldots ,g_l )^T$, and assume for simplicity that
$\displaystyle{ \sum_{j=1}^{\ell} f_j (z) g_j (z) = 0}$, $z \in
\Sigma$. Then, if $m$ solves the normalized RHP $(\Sigma ,v)$, we have
\begin{equation}\label{44.1}
F = m_{\pm} f \ \ \ \ {\rm and} \ \ \ \ G = (m_{\pm}^T )^{-1} g.
\end{equation}
The proofs of \eqref{39.1} and \eqref{41.1} proceed by expressing the relative
determinants $\frac{D_n (\omega_1 \omega_2 )}{D_n (\omega_2 )}$,
$\frac{\Delta_n(\omega_1 \omega_2 )}{\Delta_n (\omega_2 )}$ in terms of Fredholm
determinants of integrable operators $K$,
\begin{eqnarray*}
\log \det (1-K) & = & \int_0^1 {d\over dt} \log \det (1-tK) \\
 & = & -\int_0^1 \tr \left( {1\over 1-tK} K \right) dt
\end{eqnarray*}
and then using \eqref{44.1} to express $(1-tK)^{-1} K = ((1-tK)^{-1}
-1)/t$ in terms of the solution of the auxiliary RHP associated to
$tK$. We shall say more about \eqref{39.1} and \eqref{41.1} in what follows.

\section{Applications of $({\IR},v)$ and $(S^1 ,v)$: asymptotics}

In this section we consider the asymptotics of OP's, denoted (b) in
Section 2. In Section 2, the goal was to show how a variety of
identities, equations and formulae, mostly classical and well-known,
follow from a single, basic methodology in RHP's. Here the goal is
to describe new results on the asymptotics of OP's that follow from the
RH method, utilizing in particular the non-linear, non-commutative,
steepest descent method introduced in \cite{Deift-Zhou} in 1993. Although
much was known (see \cite{Szego}) about the detailed asymptotic
behavior of classical OP's, like Hermite, Laguerre, Jacobi
polynomials, etc., both on and off the contour of orthogonality,
little was known about the detailed asymptotics of OP's with respect
to general weights. The main tool that makes possible the detailed
analysis of the asymptotics of classical OP's is the existence of
integral representations for these polynomials, to which the
classical method of steepest descent can be applied (see, for
example, \cite{Szego}, Section 8.71). For general weights, one may
view the RHP's $({\IR},v(x))$ and $(S^1, v(\theta ))$ as
non-commutative analogs of these integral representations, with the
non-commutative steepest descent method now playing the role of the
classical steepest descent method.

We now describe the steepest descent method for RHP's in broad
outline: Unfortunately we do not have sufficient space in this
article to describe the method in detail. In the case of NLS (cf. \eqref{27.1} and
\eqref{27.2}), we write the solution $q(x,t)$ of the Cauchy problem for NLS
as a functional $f$, say, of the data $re^{i\theta}$,
\begin{equation}\label{48.1}
q(x,t) = f(re^{i\theta}).
\end{equation}
>From \eqref{18.1} and \eqref{27.2} we see that
\begin{equation}\label{48.2}
f(re^{i\theta} ) = \left( \int_{\IR} \mu (s;x,t) (\omega_+ +
\omega_- ) {ds\over 2\pi} \right)_{12}.
\end{equation}
Using the factorization
\begin{equation}\label{48.3}
v_{x,t} = \left( \begin{array}{ll} 1 & -re^{i\theta} \\ 0 & 1
\end{array} \right)^{-1} \left( \begin{array}{ll} 1 & 0 \\ -\bar{r}
e^{-\theta}  & 1 \end{array} \right)
\end{equation}
(cf. \eqref{16.4}), so that
\begin{equation}\label{48.4}
w_+ = \left( \begin{array}{ll} 0 & 0 \\ -\bar{r} e^{-\theta}  & 0
\end{array} \right) \ , \ \ \ \ \ w_- = \left( \begin{array}{ll} 0 &
-\bar{r} e^{-i\theta} \\ 0  & 0 \end{array} \right)
\end{equation}
we obtain
\begin{equation}\label{48.5}
q(x,t) = \left( \int_{\IR} ((I -C_{\omega}^{\IR} )^{-1} I )
(\omega_+ + \omega_- ) {ds\over 2\pi} \right)_{12}\,.
\end{equation}
For $r$ ``small", we have
$$
\begin{aligned}
q(x,t)
&= \left( \int_{\IR} ((I + \text{``small''})
(\omega_+ + \omega_- ) {ds\over 2\pi} \right)_{12}\\
&= \int_{\IR} r(s) e^{i(xs-ts^2 )} {ds \over 2\pi } + \text{``small''}
\end{aligned}
$$
indicating that the classical steepest descent method can be  applied as
$t\rightarrow\infty$. However, when $r$ is no longer ``small'', we
see from the non-linear dependence of $q(x,t)$ on $r$ in \eqref{48.5}, and
from the matrix nature of the problem, that a non-linear,
non-commutative version of the steepest descent method is required,
and this is the kind of method that was introduced in \cite{Deift-Zhou}.
In the classical steepest descent method, the integral localizes as
$t\rightarrow\infty$ to a small neighborhood of the stationary phase
point(s), $\theta^{\prime} (z_0 )= 0$, $z_0 = x/2t$ in the case of
NLS, and an explicit asymptotic formula for the solution is then
obtained by evaluating a Gaussian integral: in the fully non-linear
case (see \cite{Deift 4} \cite{DZ1}) the RHP $({\IR}, v_{x,t} )$ localizes
to a RHP in the neighborhood of the stationary phase point $z_0 =
x/2t$, and an asymptotic form for the solution
\begin{equation}\label{50.1}
q(x,t) \sim {1\over t^{1/2}} \alpha (z_0 )
  e^{i(tz_0^2 -\beta(z_0 )\log t )}
\end{equation}
is then obtained by solving this local RHP explicitly (in terms of
parabolic cylinder functions, as it turns out). The asymptotic form \eqref{50.1} was
first obtained by Zakharov and Manakov \cite{ZM}, by other means.
In situations where
there is more than one stationary phase point, for example for MKdV,
where $\theta = xz + 4tz^3$ with stationary phase points $\pm z_0 =
\pm \sqrt{-x/12t}$, the long-time behavior of solutions of MKdV (see
\cite{Deift-Zhou}) is a superposition of NLS-like contributions from
$+z_0$ and $-z_0$, as long as these points remain separated, i.e. ${-x\over t} > c >0$.
However, in the space-time
region where $-x/12t\to0$, and hence $+z_0 \rightarrow -z_0$, one is in
a non-linear ``caustic'' region which is manifested by the solution
taking the form of a self-similar oscillation, $\displaystyle{
q(x,t) \sim {1\over (3t)^{1/3} } u(x/(3t)^{1/3}) }$, where $u$
is a solution of the Painlev\'{e} II equation $u^{\prime\prime} (t) = tu +
2u^3$ (see \cite{Deift-Zhou}).

Up till this point, the RH asymptotic theory proceeded as a
non-linear analog of the classical steepest descent method in which all the
phenomena that arose could be viewed as non-linear counterparts of
phenomena that had already arisen in the linear, scalar situation. However,
with the analysis of the collisionless shock region for KdV (see
\cite{DVZ1},\cite{DZ4}), and the analysis of the asymptotic behavior of solutions of
the Painlev\'{e} II equation, it began to be clear that there were
phenomena inherent in the non-linear steepest descent method that
had no analog in the classical situation. Most importantly, it
became clear that instead of stationary phase points, one could have
``stationary phase lines'' in which case all the points on some interval
in $\IC$ contributed equally to the asymptotic behavior of the
solution of the problem. Moreover, in place of modulated linear
oscillations as in \eqref{50.1}, one would now have genuinely non-linear
oscillations described in terms of Jacobi's $sn$ and $cn$
functions, etc. A systematic extension of the steepest descent
method to allow for such ``stationary phase lines'' and genuinely
non-linear oscillations was presented by Deift, Venakides and Zhou
\cite{DVZ2} in the context of their work on the zero dispersion problem.
Soon thereafter, using the methods in \cite{DVZ2} together with recent
developments in the theory of logarithmic potentials with external
fields (see \cite{Saff-Totik}, and also \cite{DKM}), the authors in \cite{DKMVZ1}
derived so-called Plancherel-Rotach asymptotics for OP's with
measures of the form
\begin{equation}\label{53.1}
e^{-V(x)} dx,\quad V(x) = \gamma x^{2m} + \delta x^{2m-1}+\cdots
\quad\gamma > 0,
\end{equation}
and in \cite{DKMVZ2}, for measures of the form
\begin{equation}\label{53.2}
e^{-nQ(x)} dx, \quad Q(x)/\log |x| \rightarrow + \infty\quad \text{as}\quad |x| \rightarrow\infty,
\end{equation}
where $Q(x)$ is real analytic on $\IR$.
As described in \cite{DKMVZ1} one obtains as $n\rightarrow\infty$ precise
pointwise asymptotics for the OP's $P_n (z)$ for all $z\in {\IC}$,
as well as detailed asymptotics for $a_n$, $b_n$, $\gamma_n$ and the
zeros of $p_n (z)$. In the special case $e^{-n(x^4-tx^2)}dx$, Bleher and Its \cite{BI}
obtained asymptotics for the associated OP's using RH techniques and a mixture of
steepest descent/isomonodromy ideas.

In broad outline the method proceeds as follows.
For weights $e^{-V(x)}$ as above one first scales $x\to xn^{1/2m}$
so that $e^{-V(x)}\to e^{-nV_n(x)}$, where $V_n(x)=\gamma x^{2m}+\frac{\delta}{n^{1/2m}}x^{2m-1}+\cdots$.
Next, one considers the so-called equilibrium measure $d\mu_{\rm eq}$ for the
logarithmic potential problem associated with OP's (see
\cite{Saff-Totik}). By \cite{DKM}, for weights $e^{-nV_n(x)}$ or $e^{-nQ(x)}$ as
above, $d\mu_{\rm eq}$ is supported on a finite union of disjoint
intervals $\displaystyle{ \cup_{i=1}^J (a_i ,b_i ) }$, $J < \infty $
(in the case $e^{-nV_n(x)}$, $J=1$). Next one introduces the so-called
``$g$'' function, $\displaystyle{ g(z) \equiv \int_{\IR} \log (z-s)
d\mu_{\rm eq} (s) \sim \log z }$ as $z\rightarrow\infty$. Along with
$d\mu_{\rm eq}$, the logarithmic potential problem also produces a
Lagrange multiplier $\ell$, and we set $\displaystyle{
\tilde{X}^{(n)} \equiv e^{ {n\ell\over 2} \sigma_3} X^{(n)} (z)
e^{-ng(z) \sigma_3} e^{- {n\ell\over 2} \sigma_3} }$. One observes
that $\tilde{X}^{(n)}$ now solves a {\it normalized} RHP $({\IR},
\tilde{v})$ for some explicit jump matrix $\tilde{v}$. In the key
step, the RHP for $\tilde{X}^{(n)}$ is now deformed to a RHP on a
contour $\hat{\Sigma}$ of the form shown in Figure 2.
\begin{figure}
\begin{center}
\epsfig{file=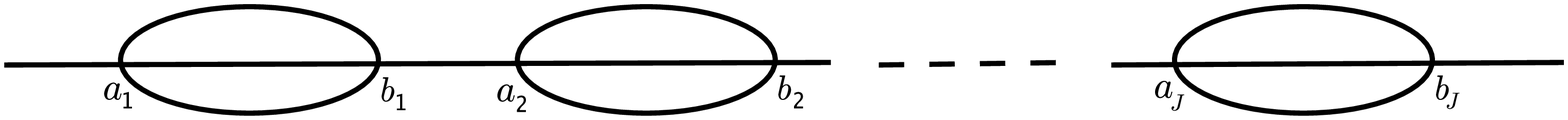,width=4.5in}
\end{center}
\caption{The contour $\hat\Sigma$}
\end{figure}
By the properties of $g(z)$, or more properly, the
properties of $d\mu_{\rm eq}$, it turns out that as
$n\rightarrow\infty$, $\hat{v}$, the jump matrix for the deformed RHP on
$\hat{\Sigma}$, converges
\begin{equation}\label{55.1}
\hat{v} (z) \rightarrow I
\end{equation}
exponentially for all $\displaystyle{ z\in \hat{\Sigma} \backslash
\cup_{i=1}^J [a_i ,b_i ] }$. Thus as $n\rightarrow\infty$, the RHP
reduces to a limiting RHP on the union of intervals $\displaystyle{
\cup_{i=1}^J [a_i ,b_i ] }$. On each of the intervals $(a_i ,b_i )$,
$\hat{v} (z)$ has the simple form $\displaystyle{ \left(
\begin{array}{ll} 0 & 1 \\ -1 & 0 \end{array} \right) }$ and this
limiting RHP can be solved explicitly in terms of the function
theory on the hyper-elliptic Riemann surface obtained by gluing
together two copies of $\displaystyle{ {\IC} \backslash
\cup_{i=1}^J (a_i ,b_i ) }$ in the standard way. However, the convergence
rate in \eqref{55.1} is not uniform, becoming slower and slower as $z$
approaches the end points $a_i ,b_i$. The natural topology for RHP's
is convergence for the coefficients of $\hat v$ in $L^p \cap L^{\infty}
(\hat\Sigma )$ (cf. \eqref{17.2}\eqref{17.5}), and the lack of uniform convergence in \eqref{55.1}
constitutes the major technical difficulty in implementing the
steepest descent method as described above. We refer the reader to
\cite{DKMVZ2, DKMVZ1} for more details.

We now consider the relative determinant formulae \eqref{39.1} and \eqref{41.1}
and their associated RHP's $({\IR} ,v_t (x))$ and $(S^1,v_t (\theta ))$ respectively.
The celebrated strong Szeg\"{o} limit
theorem, in the definitive form due to Ibragimov (see \cite{Simon 1} for many
proofs and much historical discussion) states that if $d\mu (\theta
) = e^{-V(\theta )} \displaystyle{ {d\theta \over 2\pi}}$, and
$V(\theta )$ has Fourier coefficients $\{ \hat{V}_k \}$ satisfying
$\displaystyle{ \sum_{k=1}^{\infty} k|\hat{V}_k |^2 < \infty }$,
then as $n\rightarrow\infty$
\begin{equation}\label{57.1}
\ln \Delta_n (e^{-V}) = (n+1) \hat{V}_0 + \sum_{k=1}^{\infty}
  k|\hat{V}_k |^2 + o(1) \ .
\end{equation}
In addition to the many proofs in \cite{Simon 1}, \eqref{57.1} can
also be proved, under certain additional smoothness assumptions on
$V(\theta )$, by applying the steepest descent method to the RHP's
$(S^1, v_t (\theta ))$, $0 < t < 1$. The situation is simpler than
in \cite{DKMVZ2, DKMVZ1}, but the argument in this situation is
particularly illustrative of the emergence of a ``stationary phase
line'': details are given in \cite{Deift 2}. There is also a version
of the strong Szeg\"o limit theorem for block Toeplitz determinants
(see \cite{W1}, \cite{W2}, and also \cite{Bot} for more recent results).
In the block Toeplitz case, the analog of \eqref{57.1} contains a certain Fredholm
determinant which is difficult to evaluate in elementary terms. In certain cases
the method in \cite{Deift 2} extends to the block Toeplitz case and, quite surprisingly,
the term corresponding to this Fredholm determinant is evaluated automatically
(see \cite{IJK}).

In \cite{BCW1} (see also \cite{BCW2}) the authors state the following analog
of the Szeg\"{o} strong limit theorem for the case of Hankel matrices.
Let $\omega_2 = e^{-x^2}$ and let $\omega_1 (x) > 0$ have the
property that $\omega_1 (x) \rightarrow 1$ sufficiently rapidly as
$|x| \rightarrow\infty$. Then as $n\rightarrow\infty$
\begin{equation}\label{58.1}
\ln {D_n (\omega_1 \omega_2 )\over D_n (\omega_2 )} = {
\sqrt{2(n+1)} \over \pi } \int_\IR \log \omega_1
(x)dx + {1\over 4\pi} \int_{\IR} |k| \ |\hat{f} (k)|^2 dk + o(1)
\end{equation}
where $\displaystyle{ \hat{f} (k) = {1\over \sqrt{2\pi} } \int (\log
\omega_1 (x) ) e^{-ikx} dx }$. Using \eqref{39.1}, this result can also be
proved (\cite{Deift 5}) using the steepest descent method, not only for
$\omega_2 = e^{-x^2}$, but also for more general weights, $\omega_2
= e^{-V(x)}$, $V(x) = \gamma x^{2m} + \ldots$, $\gamma > 0$, as
above.

Riemann-Hilbert techniques and the RH method are useful not only for
asymptotic evaluation, but also for estimation. For example, let
$\omega (\theta ) \in L^{\infty} (S^1 )$ be a bounded weight
on $S^1$ with Fourier coefficients $\displaystyle{ \omega_k
= \int_{-\pi}^{\pi} e^{-ik\theta} \omega (\theta ) {d\theta \over
2\pi} }$, $k \in {\ZZ}$. Let $\displaystyle{ ((T(\omega ))_{jk}
)_{j,k=0}^{\infty} = (\omega_{j-k} )_{j,k=0}^{\infty} }$ denote the
Toeplitz matrix associated with $\omega$ acting on $\displaystyle{
\ell_2^+ = \biggl\{ u = (u_0 ,u_1 ,\ldots ) : \sum_{k=0}^{\infty}
|u_k |^2 < \infty \biggr\} }$, and let $\displaystyle{ ((T_n (\omega
))_{jk} )_{j,k=0}^n = (\omega_{j-k} )_{j,k=0}^n }$ denote the
leading $(n+1) \times (n+1)$ section of $T(\omega)$. If $\omega $
is in the Wiener space $\displaystyle{ W^0 = \biggl\{ \omega :
\sum_{j=-\infty}^{\infty} |\omega_j | < \infty \biggr\} }$ with
$\omega (\theta ) > 0$, then, by a well known theorem of Krein,
$(T(\omega ))^{-1}$ exists as a bounded operator in $\ell_2^+$. The
question is the following: How closely does $(T_n (\omega ))^{-1}$ approximate
$(T(\omega ))^{-1}$ for $n$ large? Let $\nu = (\nu_k )_{k\in\ZZ}$ be a Beurling weight
(see e.g. \cite{Simon 1}): Thus $\nu_k \geq 1$,
$\nu_k = \nu_{-k}$ and $\nu_{j+k} \leq \nu_j \nu_k$ for all $j,k\in
\ZZ$. In particular, $((1+ |k| )^{\ell} )_{k \in {\ZZ}}$,
$\ell > 0$, and $(e^{\alpha |k|} )_{k\in {\ZZ}} $, $\alpha > 0$,
are Beurling weights. Define the Beurling class
$$
W_{\nu} = \biggl\{ \omega \in L^1 (S^1 ) : \sum_{j\in
{\ZZ}} \nu_j |\omega_j | < \infty \biggr\}.
$$
Clearly $W_\nu\subset W^0$ for any Beurling weight $\nu$.
Let $\omega \in W_{\nu}$ for some $\nu$ and assume in addition, for
simplicity, that the weights increase on ${\ZZ}_+$, i.e. $\nu_j
\leq \nu_k$ for $0 \leq j < k$. Then the following is true \cite{DO}: for
$n$ sufficiently large and $0 \leq j$, $k \leq n$,
\begin{equation}\label{60.1}
\left| (T_n (\omega ))_{jk}^{-1} - (T(\omega ))_{jk}^{-1} \right|
\leq c_{\nu} (\omega ) \min (\nu_{n+1-k}^{-1}, \nu_{n+1-j}^{-1} )
\end{equation}
for some constant $c_{\nu} (\omega )$.
Thus for $0\leq j,k \leq n$, $\bigl(T_n(\omega)\bigr)_{jk}^{-1}$ is a good
approximation to $\bigl(T(\omega)\bigr)_{jk}^{-1}$, apart from
the lower right corner $j\sim k\sim n$. This estimate is a
generalization of an earlier estimate due essentially to Widom (see
\cite{BS} for references and further discussion). The proof of \eqref{60.1} in
\cite{DO} uses RH techniques in an essential way closely related to
the proof of \eqref{37.2}. The paper also contains
other results for orthogonal polynomials on the unit circle,
including a new RH proof of the reverse statement in Baxter's
theorem (cf. \cite{Simon 1}). Interestingly, the Borodin-Okounkov
operator \cite{BO}, or more properly, the Borodin-Okounkov-Case-Geronimo
operator, which has emerged recently as a powerful tool in the
analysis of Toeplitz determinants, arises naturally in the analysis
in \cite{DO}.

The steepest descent method for varying weights $\omega (x) =
e^{-nQ(x)}$ in \cite{DKMVZ2} can also be applied to orthogonal
polynomials on the unit circle with varying weights $\omega (\theta
) = e^{-nQ(\theta )}$. For example in their analysis of the length $\l_n
= \l_n (\pi )$ of the longest increasing subsequence of a random
permutation $\pi$ on $n$ letters, the authors in \cite{BDJ} prove that
\begin{equation}\label{62.1}
\lim_{n\rightarrow\infty} {\rm Prob} \left( {\l_n -2\sqrt{n} \over
n^{1/6} } \leq t \right) = F_2(t)
\end{equation}
where $F_2(t)$ is the Tracy-Widom distribution function for the
largest eigenvalue of a random matrix from the Gaussian Unitary
Ensemble. The proof of \eqref{62.1} in \cite{BDJ} reduces, by a formula of
Gessel, to the analysis of the Toeplitz determinant $\Delta_{n-1}
(e^{s\cos\theta})$ where $\displaystyle{ s = (n+1) \left( 1 - {t \over
2^{1/3} ( n+1 )^{2/3} } \right) }$ as $n\rightarrow\infty$, and where
$t$ is the same as in \eqref{62.1}. As indicated above, the method of \cite{BDJ}
is modeled on the RH steepest descent method in \cite{DKMVZ2}. The same RH problem
with weight $e^{s\cos\theta}$ on $S^1$ also appears in the work of Baik and Rains \cite{BR}
in their analysis of monotone subsequences of involutions.

The steepest descent method for OP's $\{\phi_n\}$ on the unit circle can also be used to
obtain detailed information on the zeroes of the $\phi_n$'s as $n\to\infty$ (see \cite{MMS}).
In a further development, the authors of \cite{MM} have introduced an extension of the
steepest descent method to non-analytic weights, obtaining in particular new results for the
zeros of OP's on the unit circle for such weights.

Throughout this paper we have restricted our attention to measures that are smooth as in \eqref{19.1}
and \eqref{21.2}. The OP problem for general measures $d\mu$ is then analyzed (cf. \textbf{Technical Remark} above)
by approximating the measure appropriately by smooth measures $d\mu_{\epsilon}$, and then taking the limit as $\epsilon\to0$.
This approach works well for the derivation of equations, formulae, etc., but for asymptotic questions
one clearly needs a different approach. Recently remarkable connections have been discovered
(\cite{Johansson}) between various combinatorial problems - random growth models, random word problems,
tiling problems - and certain polynomials orthogonal with respect to discrete measures. The
polynomials that arise include the classical Meixner, Charlier, Krawtchouk and Hahn polynomials
(see \cite{Szego}). Related discoveries have also been made in the representation of the infinite
dimensional symmetric and unitary groups \cite{BO1}\cite{BO2}. The Meixner, Charlier and Krawtchouk
polynomials all have convenient integral representations (see \cite{Szego}) and their asymptotic
behavior can be read off using the classical method of steepest descent. This is unfortunately not
the case for the Hahn polynomials (such polynomials are needed in particular to describe the tiling
of hexagons by rhombi). It turns out, however, that discrete OP problems can be rephrased in terms
of a discrete RHP, which is an analogue of the continuous case, and which was introduced by Borodin,
along with a theory of discrete integrable operators, in \cite{B1}. In a significant further development
of the nonlinear steepest descent method, the authors in \cite{BKMM} extended the method to a wide class
of discrete RHP's which includes the discrete RHP for the Hahn polynomials (as well as the other three
discrete OP systems mentioned above). The relevant limit here is when the order of the OP's $p_n$
becomes large and simultaneously the spacing between the points in the measures goes to zero
at a prescribed rate (see \cite{BKMM}). In this way the authors are able to analyze the Hahn
polynomials asymptotically, proving \textit{en route} a conjecture of Johansson in \cite{Johansson}
that for hexagonal tiling the so-called ``arctic circle'' of \cite{CLP} exhibits Tracy-Widom fluctuations
as in \eqref{62.1} above. In \cite{BO2} the authors also consider an asymptotic problem for Hahn polynomials using
a discrete RHP, but the relevant limit is different from that in \cite{BKMM}.

Many researchers are currently involved in the application
of RH techniques to the theory of OP's. In addition to those mentioned above, the list
includes Chen, Claeys, Kapaev, Kitaev, Kuijlaars, van Assche and Vanlessen, amongst many others.
Because of space limitations, however, we unfortunately cannot describe their work in any detail,
and we must refer the reader to the literature.

\section{Related areas}

In this final section we will describe, very briefly, various areas related to OP's in which the RH
method plays a role.

We first consider random matrix theory (RMT), which has been a major source of questions and challenges to
OP theorists for over 40 years (see e.g. \cite{Mehta} and \cite{Deift 1}). The situation is as follows.
A Unitary Ensemble (UE) is an ensemble of $N\times N$ Hermitian matrices $\{M=M^*\}$ with
probability distribution
\begin{equation}\label{69.1}
P_N(M)\,dM=\frac{1}{Z_N} e^{-\tr W(M)} dM
\end{equation}
where
\begin{itemize}
 \item $dM$ denotes Lebesgue measure on the algebraically independent elements of $M$.
 \item $W(x)$ is a real-valued function that goes to $+\infty$ as $|x|\to\infty$. The case
       $W(x)=x^2$ gives rise to the Gaussian Unitary Ensemble (GUE).
 \item $Z_N$ is a normalization coefficient.
\end{itemize}

``Unitary'' refers to the fact that the distribution \eqref{69.1} is invariant
under unitary conjugation, $M\to UMU^*$, $U$ unitary. The Universality Conjecture
for UE's (see \cite{Mehta} and \cite{Deift 1}) states, in particular, the following:
Given $W$, if $J_N=c_N+s_N(-t,t)$ is a suitably centered and scaled interval in $\IR$,
then as $N\to\infty$, $P(J_N)=\text{Prob}(M : M \,\text{has no eigenvalues in}\, J_N)$
converges to a \textit{universal limit} independent of $W$,
\begin{equation}\label{70.1}
\lim_{N\to\infty} P(J_N)=\det{(1-S_t)}
\end{equation}
where $S_t$ is the trace class operator with kernel $S_t(x,y)=\frac{\sin\pi(x-y)}{\pi(x-y)}$
acting in $L^2(-t,t)$. The specific form of the weight $e^{-W(x)}\,dx$ is reflected only
in the precise values of $c_N$ and $s_N$. OP's enter the picture because of the celebrated result of
Gaudin and Mehta (see \cite{Mehta}) that if $B\subset\IR$ is a Borel set, then
\begin{equation}\label{71.1}
\text{Prob}(M : M\,\text{has no eigenvalues in}\, B)=\det(1-K_{N,B})
\end{equation}
where $K_{N,B}$ is the finite rank operator with kernel
\begin{equation}\label{71.2}
K_N(x,y)=\sum_{j=0}^{N-1} p_j(x)p_j(y) e^{-\frac12 W(x)} e^{-\frac12 W(y)}
\end{equation}
acting on $L^2(B)$, and $\{p_j\}_{j\geq0}$ are the orthonormal
polynomials \eqref{2.3} with respect to the weight $e^{-W(x)}\, dx$. Hence the question
of proving universality as in \eqref{70.1} becomes a question of deriving the appropriate
asymptotics for OP's, and this is the main scientific content of \cite{DKMVZ1}, \cite{DKMVZ2},
\cite{Deift 1} and \cite{BI}. Of course, if the weight $e^{-W(x)}\,dx$ is classical, e.g.
$W(x)=x^2$, and the asymptotics of the associated polynomials $\{p_j\}_{j\geq0}$ can
be derived from an integral representation, then universality for these ensembles
can be proved without recourse to the RH steepest descent method, and this has been done
by various authors (see \cite{DKMVZ1}, \cite{DKMVZ2} for references to the literature).

Orthogonal ensembles (OE's) of $N\times N$ real symmetric matrices
$\{M=\bar M=M^T\}$ and Symplectic Ensembles (SE's) of $2N\times2N$ Hermitian self-dual
matrices $\{M=M^*, JMJ^T=M^T\}$, where $J=\diag(\tau,\ldots,\tau)$,
$\tau=\left(
\begin{matrix}0 &1\\ -1&0
\end{matrix}\right)$,
equipped with invariant weights analogous to \eqref{69.1}, are more difficult to
analyze. Firstly, in the place of determinantal expressions as in \eqref{71.1},
one obtains Pfaffians (see \cite{Mehta} for classical ensembles, \cite{TW} for the
general case)
\begin{equation}\label{73.1}
\text{Prob}(M : M\,\text{has no eigenvalues in}\, B)=\bigl(\det(1-\hat K_{N,B})\bigr)^{1/2},
\end{equation}
and, moreover, the operators $\hat K_{N,B}$ are now $2\times2$ matrix operators with kernels
$\bigl(\hat K_{N,ij}(x,y)\bigr)_{1\leq i,j\leq2}$, $x,y\in B$. In contrast to \eqref{71.2},
these kernels are most naturally expressed in terms of certain skew-orthogonal polynomials
(see \cite{Mehta}), but for general weights $e^{-W(x)}\,dx$ the asymptotic behavior of such
polynomials is not known. However Widom \cite{W} has shown that if $W'/W$ is rational, then
$\bigl(\hat K_{N,ij}(x,y)\bigr)$ can be expressed conveniently in terms of the orthonormal
polynomials $\{p_j\}_{j\geq0}$ with respect to the weight $e^{-W(x)}\,dx$, so again, as in the unitary case,
the question of universality of OE's and SE's becomes a question of analyzing the asymptotic behavior of OP's.
The expressions for $\bigl(\hat K_{N,ij}(x,y)\bigr)_{1\leq i,j\leq2}$ are now more cumbersome than \eqref{71.2} and significant
new technical issues arise, but nevertheless, using the asymptotic analysis in \cite{DKMVZ1} as a basic
ingredient, it is indeed possible to use Widom's formulae in \cite{W} to prove universality for OE's and SE's
with weights of the form $e^{-V(x)}\,dx$, $V(x)=\gamma x^{2m}+\cdots$, $\gamma>0$. This is the content of
\cite{DG1} and \cite{DG2}.

Biorthogonal polynomials $\pi_k(x)=x^k+\cdots$, $\sigma_j(y)=y^j+\cdots$, $k,j\geq0$,
\begin{equation}\label{75.1}
\int_\IR \int_\IR \pi_k(x)\sigma_j(y) e^{-V(x)-W(y)+2\tau xy}\,dxdy=0\quad\text{if}\,j\neq k,
\end{equation}
arise in the analysis of the theory of coupled random matrices. Here $V(x)$ and $W(y)$
grow sufficiently rapidly as $|x|,|y|\to\infty$, and $\tau\neq0$. Various RH problems have been
proposed to analyze these polynomials (see, in particular, \cite{BEH}, \cite{Kap}, \cite{KM}
and the references therein), but the analysis of the asymptotic behavior of these RHP's is still
at a preliminary stage.

For $m\geq 2$, let $n=(n_1,n_2,\ldots,n_m)$ be a vector of non-negative integers, and let
$\omega_1(x)\geq0,\ldots, \omega_m(x)\geq0$ be weights on $\IR$ with finite moments. Let
$|n|=n_1+\cdots+n_m$. Multiple orthogonal polynomials (see \cite{Apt}) of type I are polynomials
$A_n^{(k)}$ for $k=1,2,\ldots,m$, $\text{deg}\, A_n^{(k)}\leq n_k-1$ such that the function
$$
H_n(x)=\sum_{k=1}^m A_n^{(k)}(x)\omega_k(x)
$$
satisfies
\begin{equation}\label{76.1}
\int_\IR x^j H_n(x)\, dx=
\left\{
  \begin{array}{ll}
    0, & \text{for}\, j=0,\ldots,|n|-2; \\
    1, & \text{for}\, j=|n|-1.
  \end{array}
\right.
\end{equation}
Multiple orthogonal polynomials $L_n(x)$ of type II are monic polynomials of degree $|n|$
satisfying
\begin{equation}\label{77.1}
\int_\IR L_n(x) x^k\omega_j(x)\, dx=0\quad\text{for}\,\,k=0,\ldots,n_j-1,\,j=1,\ldots,m.
\end{equation}
Multiple orthogonal polynomials were first introduced by Hermite in his proof of the
transcendence of $e$. In 2000, van Assche, Geronimo and Kuijlaars \cite{vAGK} showed that
multiple orthogonal polynomial problems of types I and II could be rephrased
as RHP's analogous to the RHP of Fokas, Its and Kitaev for ordinary OP's, and they used
these RHP's to derive various properties and relations for the multiple OP's. In the last
year or two significant progress has been made in extending and applying the steepest descent method to
RHP's which arise from multiple OP's in special cases. We mention, in particular, \cite{BK}, \cite{ABK}
and \cite{KVaW}, \cite{KSVaW} and the references therein: In the first two papers the authors
consider a random matrix ensemble $P_N(M)\,dM=\frac{1}{Z_N} e^{-N\tr(\frac12 M^2-AM)}\,dM$,
with external source $A$, first analyzed by Pastur, Br\'ezin-Hikami, and later by Zinn-Justin.
Under certain conditions on $A$, they show that the ensemble can be analyzed as
$N\to\infty$ in terms of a $3\times3$ RHP to which an extension of the nonlinear steepest
descent method can be applied: A new phenomenon now occurs in the analysis, which the authors
term a ``global opening of lenses'' (see \cite{ABK}). In the second two papers the authors
analyze type I and type II Hermite-Pad\'e approximations to the exponential function, which
they are again able to control by applying an extension of the steepest descent to a
$3\times3$ RHP.

Riemann-Hilbert techniques can also be used to analyze the asymptotics of so-called
\textit{orthogonal Laurent polynomials}. Such polynomials arise in the following way.
Let $V(x)$ be a real-analytic function on $\IR\setminus \{0\}$ with the property
$$
\lim_{|x|\to\infty} \frac{V(x)}{\ln|x|}=\lim_{|x|\to0} \frac{V(x)}{\ln(|x|^{-1})}=+\infty.
$$
Orthogonalization of the ordered basis $\{1,z^{-1},z,z^{-2},z^2,\ldots\}$ with
respect to the pairing $(f,g)\mapsto\int_\IR f(s)g(s) e^{-NV(s)}\,ds$ leads to the even degree
and odd degree orthonormal Laurent polynomials $\{\phi_m\}_{m\geq0}$:
$\phi_{2n}(z)=\xi^{(2n)}_{-n}z^{-n}+\cdots+\xi^{(2n)}_{n}z^{n}$, $\xi^{(2n)}_n>0$,
$\phi_{2n+1}(z)=\xi^{(2n+1)}_{-n-1} z^{-n-1}+\cdots+\xi^{(2n+1)}_{n}z^{n}$,
$\xi^{(2n+1)}_{-n-1}>0$. Recently, McLaughlin, Vartanian and Zhou (see \cite{MVZ}
and the references therein) have used RHP-steepest descent methods to analyze the
asymptotic behavior of the Laurent polynomials $\phi_{2n}(z), \phi_{2n+1}(z)$ and
their associated norming constants $\xi^{(2n)}_n, \xi^{(2n+1)}_{-n-1}$
in the limit as $N\to\infty$, $N/n\to1$. The work of McLaughlin et al. involves
significant extensions of the steepest descent method: Such extensions are needed in order to overcome the
new difficulties introduced into the problem by the singularity of the potential $V(x)$
at $x=0$.

Finally, there are problems in which the asymptotic behavior of the system at hand is
described by OP's. This happens, in particular, in the case of the so-called Toda
rarefaction problem (see \cite{DKKZ}). Here one considers the initial-boundary value problem
for the Toda lattice
\begin{equation}\label{79.1}
\ddot{x}_n=e^{x_{n-1}-x_n}-e^{x_{n}-x_{n+1}},\quad n\geq1
\end{equation}
where for some $\alpha>0$
\begin{equation}\label{79.2}
\left\{
  \begin{array}{ll}
    x_n(0)=\alpha n, & n\geq1; \\
    \dot x_n(0)=0, & n\geq1,
  \end{array}
\right.
\end{equation}
and the driving particle moves with a fixed velocity $2a$
\begin{equation}\label{79.3}
x_0(t)=2at,\qquad t\geq0.
\end{equation}
Making the change of variables $x_n\to \alpha n+y_n$ one sees that, apart from rescaling time,
one can always assume without loss of generality that $\alpha=0$ in \eqref{79.2}. One may think of
\eqref{79.1}--\eqref{79.3} as a cylinder of particles $\{x_n\}_{n\geq1}$ driven by a piston $x_0$.
If $a>0$, one has the (Toda) shock problem (\cite{VDO}) and if $a<0$ one has the (Toda) rarefaction
problem. In the rarefaction problem, if $|a|$ is sufficiently large ($|a|>1$ turns out to be the
critical region) one expects that the piston will separate from the ``gas'' $\{x_n\}_{n\geq1}$
and cavitation will occur. This is indeed what happens: if $a<-1$, the authors in \cite{DKKZ}
show, using the RH steepest descent method, that as $t\to\infty$, the solution of the Toda lattice
splits into two parts, I+II. Part I models the cavitation and Part II is an exponentially
decreasing error term. Quite remarkably, Part I is constructed from the solution of an
associated OP problem, which turns out to be the Fokas, Its, Kitaev RHP in disguise. We refer the reader
to \cite{DKKZ} for details.


\end{document}